# HAZARD RATE ESTIMATION FOR LOCATION-SCALE DISTRIBUTIONS UNDER COMPLETE AND CENSORED DATA

**Baris Surucu**
**Department of Statistics**
**Middle East Technical University**
**06800 Ankara / Turkey**

## 1. Introduction

In reliability and life testing studies, the topic of estimating hazard rate has received great attention in recent years since an estimate of hazard rate is a quite useful tool for making decisions. Some works have included nonparametric approaches while some have considered parametric structural models for complete as well as censored data sets; see Meeker et al. (1992), Antoniadis and Grégoire (1999), Rai and Singh (2003), Bezandry et al. (2005), Brunel and Comte (2008), and Mahapatra et al. (2012). Depending on the shapes of the hazard rate, efficiencies differ markedly across proposed estimators. This situation is remarkable especially when different estimation techniques are utilized for unknown parameters of underlying distributions in parametric approaches. That is, estimated hazard rate (and also reliability) at a specific time point $t$ as functions of these estimators leads to inconsistent coverage probabilities as distributional convergence of hazard rate estimator may not be at the desired rate for certain sample sizes.

This manuscript focuses on the estimation of monotone hazard rate when the distribution of concern is of location-scale type. A very simple and efficient approximation of hazard rate for a complete sample is introduced as a function of order statistics, which allows a fast convergence to the asymptotic distribution while achieving highly accurate coverage probabilities for confidence intervals (CI). Efficiency of the method is demonstrated on simulation studies by considering a number of location-scale distributions having various hazard shapes. We also consider the case of censored samples by incorporating the predicted future order statistics into the picture. The remainder of the paper is organized as follows: In Chapter 2, we discuss the estimation of unknown parameters for location-scale families and their efficiencies. In Chapter 3, we propose a linear approximation of the hazard rate function to estimate the hazard rate at time $t$ for different underlying location-scale distributions. We discuss the asymptotic properties of the hazard rate estimator and present the efficiencies of the new method by comparing with those of classical approach in Chapter 4. For the censored samples, necessity of predicting future order statistics in the novel method is explained in Chapter 5. Moreover, a second approach based on expected spacings is introduced to achieve the same goal. In the same chapter, some graphical displays are provided to illuminate the discussion on boundaries of CI's for hazard rate before concluding.

## 2. Parameter Estimation

Let us have a location-scale distribution having the probability density function

$$f_X(x;\mu,\sigma) = \frac{1}{\sigma} f\left(\frac{x-\mu}{\sigma}\right) \qquad (1)$$

where $\mu$ and $\sigma$ are location and scale parameters, respectively. For (1), the hazard rate at time point $t$ can be written as

$$h\left(\frac{t-\mu}{\sigma}\right) = f\left(\frac{t-\mu}{\sigma}\right) \Big/ \left\{1 - F\left(\frac{t-\mu}{\sigma}\right)\right\} \qquad (2)$$

where $f(z)$ and $F(z)$ are probability density function and distribution function of the standardized random variable $z = \frac{x-\mu}{\sigma}$, respectively. To estimate the hazard rate, we firstly need to estimate the unknown location and scale parameters. Throughout the work, we will consider two types of estimation techniques; least squares (LS) estimation and modified maximum likelihood (MML) estimation of Tiku (1967). LS estimators are well-known and we do not give the details about them in this manuscript. To briefly introduce MML, write the log-likelihood function of a censored sample ($n$-$r$ largest observations censored) as

$$\ln L = -r \ln \sigma + \sum_{i=1}^{r} \ln f\left(\frac{x_i - \mu}{\sigma}\right) + (n-r) \ln\left[1 - F\left(\frac{x_{(r)} - \mu}{\sigma}\right)\right]. \qquad (3)$$

Realize that when $r=n$, we have the complete sample of size $n$. Then, the partial derivatives with respect to the unknown location and scale parameters can be written as

$$\frac{\partial \ln L}{\partial \mu} = -\frac{1}{\sigma} \sum_{i=1}^{r} \frac{\frac{\partial}{\partial \mu} f\left(\frac{x_i - \mu}{\sigma}\right)}{f\left(\frac{x_i - \mu}{\sigma}\right)} - \frac{(n-r)}{\sigma} \frac{f\left(\frac{x_{(r)} - \mu}{\sigma}\right)}{1 - F\left(\frac{x_{(r)} - \mu}{\sigma}\right)} \qquad (4)$$

and

$$\frac{\partial \ln L}{\partial \sigma} = -\frac{r}{\sigma} - \frac{1}{\sigma^2} \sum_{i=1}^{r} \frac{\frac{\partial}{\partial \sigma} f\left(\frac{x_i - \mu}{\sigma}\right)}{f\left(\frac{x_i - \mu}{\sigma}\right)} + \frac{(n-r)}{\sigma^2} \frac{f\left(\frac{x_{(r)} - \mu}{\sigma}\right)}{1 - F\left(\frac{x_{(r)} - \mu}{\sigma}\right)}. \qquad (5)$$

Note also that sample observations $x_i$ in (4) and (5) can be replaced by ordered statistics $x_{(i)}$ as the complete sums are invariant to ordering. Then, define $g_1(z_{(i)}) = \dfrac{\frac{\partial}{\partial \mu} f(z_{(i)})}{f(z_{(i)})}$, $g_2(z_{(i)}) = \dfrac{\frac{\partial}{\partial \sigma} f(z_{(i)})}{f(z_{(i)})}$ and $g_3(z_{(r)}) = \dfrac{f(z_{(r)})}{1 - F(z_{(r)})}$, where $z_{(i)}$ is the $i$th standardized order statistic. Depending on the structural form of the pdf of the random variable X, the functions $g_1, g_2$ and $g_3$ are mostly nonlinear and the basic cause of having inexplicit solutions for the likelihood functions. To overcome this problem, Tiku (1967) proposed a modification of the maximum likelihood approach, which leaded to very efficient explicit solutions for the unknown location and scale parameters of underlying distributions. In fact, the MML estimators are robust and asymptotically fully efficient; see also Tiku and Suresh (1992).

Let one of these nonlinear functions be $g(z_{(i)})$. Then, a Taylor expansion around the expected value of the $i$th standardized order statistic ($m_i$) gives us the linear approximation

$$g(z_{(i)}) \cong \alpha + \beta\, z_{(i)} \tag{6}$$

where $\alpha$ and $\beta$ are the corresponding constants obtained from the method; see Tiku (1967) and Tiku and Suresh (1992). By replacing all nonlinear functions with their approximate linear forms, we obtain explicit solutions of the equations (4) and (5), which lead to Tiku's MML estimators. In fact, the MML estimators are asymptotically equivalent to ML estimators and are very highly efficient. Specific forms of the estimators and their distributional properties will be examined in the next sections while working with different distributions.

## 3. Estimating Hazard Rate

The hazard rate at time $t$ for a location-scale distribution can be defined as

$$h(\delta) = f(\delta)/(1 - F(\delta)) \tag{7}$$

where $\delta = \dfrac{t - \mu}{\sigma}$. Since the parameters are unknown in the real life, we first estimate them by using a proper estimation technique and plug them in the hazard function (7) to estimate the hazard rate. However, this brings about some distributional problems for the estimator $\hat{h}(\delta)$. That is, the asymptotic distribution of $\hat{h}(\delta)$ cannot be obtained easily and, therefore, it is quite difficult to construct an CI with a relevant coverage probability.

An approximation technique similar to that discussed in Section 2 helps us to overcome these distributional problems. The new linear approximation method does not only solve the distributional struggle, but also provides a very simple and practical structural form with very high efficiency. To explain the method, consider the following linear approximation.

$$f\left(\frac{t-\hat{\mu}}{\hat{\sigma}}\right) \Big/ \left\{1 - F\left(\frac{t-\hat{\mu}}{\hat{\sigma}}\right)\right\} \cong \alpha^* + \beta^* \hat{\delta} \tag{8}$$

where $\hat{\delta} = \left(\frac{t-\hat{\mu}}{\hat{\sigma}}\right)$. Differently from Section 2, $\delta$ is not a standardized order statistic which can be expanded around its expectation by Taylor series approach. On the other hand, $t$ is always between two consecutive order statistics (for example $x_{(k)}$ and $x_{(k+1)}$). Therefore, the expected values of the standardized order statistics $z_{(k)}$ and $z_{(k+1)}$ can be considered as the boundaries of a new linearization for $\hat{h}(\delta)$. Then, the estimated hazard rate at point $t$ can be written as

$$\hat{h}(\delta) = \alpha^* + \beta^* \hat{\delta} \tag{9}$$

where

$$\alpha^* = h(m_k) - \beta^* m_k$$

$$\beta^* = \frac{h(m_{k+1}) - h(m_k)}{m_{k+1} - m_k}; \tag{10}$$

$m_k$ being the expected value of the $k$th standardized order statistic $z_{(k)}$.

**Asymptotic Distribution of the Hazard Rate Estimator:** Due to the properties of MML estimators, the asymptotic distribution of $\sqrt{m}\left(\frac{\hat{\mu}-t}{\hat{\sigma}}\right)$ is normal with mean $\sqrt{m}\left(\frac{\mu-t}{\sigma}\right)$ and variance 1, where $V(\hat{\mu}) = \sigma^2/m$; m is a constant taking different values for different distributions. Accordingly, $h(\hat{\delta}) = \alpha^* - \beta^*\left(\frac{\hat{\mu}-t}{\hat{\sigma}}\right)$ is asymptotically normally distributed with mean $\alpha^* - \beta^*\left(\frac{\mu-t}{\sigma}\right)$ and variance $\beta^{*2}/m$. From here, it is very easy to construct $(1-\alpha)\%$ CI's for both MML and LS estimators. For monotone hazard structures, the CI's are given as follows:

MML: $P\left\{\alpha^* - \beta^*\left(\dfrac{\hat{\mu}-t}{\hat{\sigma}}\right) - z_{\alpha/2}\dfrac{\beta^*}{\sqrt{m}} \leq h(\delta) \leq \alpha^* - \beta^*\left(\dfrac{\hat{\mu}-t}{\hat{\sigma}}\right) + z_{\alpha/2}\dfrac{\beta^*}{\sqrt{m}}\right\} = 1 - \alpha$, $\beta > 0$ (11)

LSE: $P\left\{\alpha^* - \beta^*\left(\dfrac{\bar{x}-t}{s}\right) - z_{\alpha/2}\dfrac{\beta^*}{\sqrt{n}} \leq h(\delta) \leq \alpha^* - \beta^*\left(\dfrac{\bar{x}-t}{s}\right) + z_{\alpha/2}\dfrac{\beta^*}{\sqrt{n}}\right\} = 1 - \alpha$, $\beta > 0$; (12)

$\bar{x}$ and $s$ being the least squares estimators of $\mu$ and $\sigma$, respectively.

**Remark:** In general, the MML-interval will on the average be shorter than the LSE-interval since m is greater than n.

## 4. Hazard Rate Estimation for Several Distributions

In this section, we will work on the hazard rate estimation process for short-tailed and long-tailed symmetric distributions which have monotone hazard rates. Performances of the estimators will be analyzed through simulation studies. For the comparison purpose, two basic types of estimation procedures will be considered. The first one is the hazard rate estimator obtained through replacing unknown parameters with their estimators (MML or LS) in the hazard rate function. The second one is due to the new linear approximation technique which utilizes two consecutive order statistics. The estimators obtained from the first approach will be named as $HR^1_{MML}$ and $HR^1_{LS}$ whereas those obtained from the second approach will be named as $HR^2_{MML}$ and $HR^2_{LS}$.

### 4.1. Short-Tailed Symmetric Distribution

Short-tailed symmetric (STS) distribution introduced by Tiku and Vaughan (1999) has the pdf

$$f(y) = C\left\{1 + \dfrac{\lambda}{2r}z^2\right\}^r \dfrac{1}{\sqrt{2\pi}}\exp\{-z^2/2\}, \quad -\infty < z < \infty,\qquad(13)$$

where $z = (y-\mu)/\sigma$; $r$ being an integer, $\lambda = \dfrac{r}{r-d}$, $r > d$ and

$$C = \dfrac{1}{\sum_{j=0}^{r}\binom{r}{j}\left(\dfrac{\lambda}{2r}\right)^j \dfrac{(2j)!}{2^j(j!)}}.\qquad(14)$$

The STS distribution is unimodal for $d<0$. For more details about the STS family, see Tiku and Vaughan (1999).

The MML estimators of a complete sample of size $n$ for the STS distribution are given as (Tiku and Akkaya, 2004, p.68)

$$\hat{\mu} = \sum_{i=1}^{n} \beta_{1i} y_{(i)} / m, \quad \left( m = \sum_{i=1}^{n} \beta_{1i} \right), \tag{15}$$

$$\hat{\sigma} = \left\{ -B + \sqrt{B^2 + 4nC} \right\} / 2\sqrt{n(n-1)}, \tag{16}$$

where

$$B = \lambda \sum_{i=1}^{n} \alpha_{1i} y_{(i)}, \quad C = \frac{2p}{k} \left\{ \sum_{i=1}^{n} \beta_{1i} y_{(i)}^2 - m\hat{\mu}^2 \right\},$$

$$\alpha_{1i} = \frac{(\lambda/r)t_i^3}{\left[1 + (\lambda/2r)t_i^2\right]^2} \quad \text{and} \quad \beta_{1i} = 1 - \lambda \frac{1 - (\lambda/2r)t_i^2}{\left[1 + (\lambda/2r)t_i^2\right]^2}. \tag{17}$$

When $\lambda \leq 1$, all $\beta_{1i}$ are positive and the equations in (17) can be used. When $\lambda > 1$, however, some $\beta_{1i}$ coefficients take negative values. To overcome the inconvenience in the estimation step, $\alpha_{1i}$ and $\beta_{1i}$ are replaced by $\alpha_{2i}$ and $\beta_{2i}$;

$$\alpha_{2i} = \frac{(\lambda/r)t_i^3 + (1 - 1/\lambda)t_i}{\left\{1 + (\lambda/2r)t_i^2\right\}^2} \quad \text{and} \quad \beta_{2i} = 1 - \lambda \frac{(1/\lambda) - (\lambda/2r)t_i^2}{\left\{1 + (\lambda/2r)t_i^2\right\}^2}; \tag{18}$$

The $(1-\alpha)\%$ MML CI for the hazard rate at time $t$ can be constructed as in the equation (10) with

$$m = n\left\{1 - \frac{2}{h}\left[\left(1 - \frac{1}{2h}\right) \Big/ \left(1 + \frac{1}{h} + \frac{3}{4h^2}\right)\right]\right\}; \quad h = 2 - d, \quad d < 2; \tag{19}$$

see also Akkaya and Tiku (2008) for more details. Since the exact expected values of the standardized order statistics are not available for the STS distribution, we use their simulated values for the linear approximation method. Calculations for $\alpha^*$ and $\beta^*$ are straightforward from (10).

Least squares estimators of $\mu$ and $\sigma$ can be obtained as

$$\bar{x} = \frac{\sum_{i=1}^{n} x_i}{n} \quad \text{and} \quad s = \sqrt{\frac{\sum_{i=1}^{n}(x_i - \bar{x})^2}{(n-1)\mu_2}}, \tag{20}$$

where

$$\mu_2 = \frac{\sum_{j=0}^{r} \binom{r}{j}\left(\frac{\lambda}{2r}\right)^j \frac{\{2(j+1)\}!}{2^{j+1}(j+1)!}}{\sum_{j=0}^{r} \binom{r}{j}\left(\frac{\lambda}{2r}\right)^j \frac{(2j)!}{2^j (j)!}} \tag{21}$$

is the variance of the short-tailed symmetric family. It is again very easy to construct the $(1-\alpha)\%$ LS CI for the hazard rate at time $t$ by simply following the rules given in Section 3.

In Table 1, one can see the simulation results for the efficiencies of the hazard rate estimators for the STS distribution. For conciseness, we only give the results for n=20. The efficiency behaviors for other sample sizes and other parameter values are almost the same. Estimators based on the MML are more efficient than those based on the LS. $HR^1_{MML}$ and $HR^2_{MML}$ have almost the same efficiency for the all quantiles.

In Table 2, the coverage probabilities for the CI's are given for $HR^2_{MML}$ and $HR^2_{LS}$. It is very clear from the simulation results that the linear approximation utilizing the MML estimators is performing very much better than that using the LS estimators. The expected $(1-\alpha)=0.95$ level is achieved in the mid-quantiles of the distribution. Towards the tails of the STS distribution, the coverage probability decreases to 0.85. On the other hand, the coverage probability for $HR^2_{LS}$ falls below the level of 0.70 in the lower and upper tails of the distribution. One can also see the Figure 1 and Figure 2 which show the CI boundaries of the MML-intervals for the same distribution.

**Table 1:** Means and Variances of Hazard Rate Estimators for the Quantiles of the STS Distribution (Q) ($n = 20$, $\mu = 0, \sigma = 1, r = 2, d = 0$)

$\hat{\mu} = 0.0068$, $\hat{\sigma} = 1.006$, $\bar{x} = 0.0064$, $s = 0.984$

| | | Means for n=20 | | | | Variances for n=20 | | | |
|---|---|---|---|---|---|---|---|---|---|
| Q | Exact | $HR^1_{MML}$ | $HR^1_{LS}$ | $HR^2_{MML}$ | $HR^2_{LS}$ | $HR^1_{MML}$ | $HR^1_{LS}$ | $HR^2_{MML}$ | $HR^2_{LS}$ |
| 0.07 | 0.1214 | 0.1238 | 0.1184 | 0.1212 | 0.1134 | 0.0032 | 0.0034 | 0.0037 | 0.0044 |
| 0.12 | 0.1798 | 0.1778 | 0.1721 | 0.1773 | 0.1706 | 0.0036 | 0.0041 | 0.0038 | 0.0045 |
| 0.17 | 0.2271 | 0.2228 | 0.2174 | 0.2228 | 0.2170 | 0.0035 | 0.0041 | 0.0036 | 0.0043 |
| 0.21 | 0.2671 | 0.2622 | 0.2574 | 0.2621 | 0.2572 | 0.0032 | 0.0038 | 0.0031 | 0.0039 |
| 0.26 | 0.3022 | 0.2979 | 0.2938 | 0.2970 | 0.2932 | 0.0028 | 0.0035 | 0.0027 | 0.0034 |
| 0.31 | 0.3341 | 0.3314 | 0.3282 | 0.3292 | 0.3269 | 0.0025 | 0.0032 | 0.0024 | 0.0031 |
| 0.36 | 0.3649 | 0.3644 | 0.3622 | 0.3603 | 0.3598 | 0.0025 | 0.0032 | 0.0024 | 0.0030 |
| 0.40 | 0.3966 | 0.3986 | 0.3977 | 0.3925 | 0.3943 | 0.0029 | 0.0037 | 0.0028 | 0.0034 |
| 0.45 | 0.4317 | 0.4366 | 0.4371 | 0.4282 | 0.4327 | 0.0037 | 0.0048 | 0.0036 | 0.0045 |
| 0.50 | 0.4728 | 0.4806 | 0.4828 | 0.4705 | 0.4776 | 0.0052 | 0.0068 | 0.0051 | 0.0066 |
| 0.55 | 0.5224 | 0.5331 | 0.5374 | 0.5216 | 0.5316 | 0.0076 | 0.0102 | 0.0076 | 0.0100 |
| 0.60 | 0.5830 | 0.5967 | 0.6035 | 0.5840 | 0.5970 | 0.0114 | 0.0155 | 0.0115 | 0.0154 |
| 0.64 | 0.6567 | 0.6733 | 0.6832 | 0.6603 | 0.6762 | 0.0170 | 0.0231 | 0.0173 | 0.0232 |
| 0.69 | 0.7456 | 0.7652 | 0.7786 | 0.7525 | 0.7713 | 0.0248 | 0.0337 | 0.0256 | 0.0340 |
| 0.74 | 0.8526 | 0.8750 | 0.8927 | 0.8633 | 0.8851 | 0.0355 | 0.0478 | 0.0367 | 0.0481 |
| 0.79 | 0.9819 | 1.0070 | 1.0299 | 0.9974 | 1.0223 | 0.0499 | 0.0662 | 0.0514 | 0.0663 |
| 0.83 | 1.1415 | 1.1692 | 1.1981 | 1.1618 | 1.1906 | 0.0692 | 0.0905 | 0.0706 | 0.0897 |
| 0.88 | 1.3477 | 1.3780 | 1.4151 | 1.3723 | 1.4065 | 0.0960 | 0.1320 | 0.0964 | 0.1204 |
| 0.93 | 1.6470 | 1.6801 | 1.7264 | 1.6739 | 1.7166 | 0.1376 | 0.1709 | 0.1345 | 0.1646 |

**Table 2:** Coverage Probabilities $(1-\alpha)$ of the CI's for the Hazard Rate Estimators for the STS Distribution ($n = 20$, $\alpha = 0.05$, $\mu = 0, \sigma = 1, r = 2, d = 0$)

$\hat{\mu} = 0.0068$, $\hat{\sigma} = 1.006$, $\bar{x} = 0.0064$, $s = 0.984$

| | Coverage Probabilities | | | | |
|---|---|---|---|---|---|
| Q | $HR^2_{MML}$ | $HR^2_{LS}$ | Q | $HR^2_{MML}$ | $HR^2_{LS}$ |
| 0.07 | 0.843 | 0.686 | 0.55 | 0.947 | 0.814 |
| 0.12 | 0.874 | 0.719 | 0.60 | 0.945 | 0.811 |
| 0.17 | 0.891 | 0.746 | 0.64 | 0.938 | 0.806 |
| 0.21 | 0.910 | 0.769 | 0.69 | 0.930 | 0.798 |
| 0.26 | 0.919 | 0.782 | 0.74 | 0.921 | 0.787 |
| 0.31 | 0.931 | 0.796 | 0.79 | 0.910 | 0.772 |
| 0.36 | 0.939 | 0.804 | 0.83 | 0.896 | 0.754 |
| 0.40 | 0.943 | 0.809 | 0.88 | 0.876 | 0.730 |
| 0.45 | 0.947 | 0.812 | 0.93 | 0.847 | 0.696 |
| 0.50 | 0.949 | 0.813 | | | |

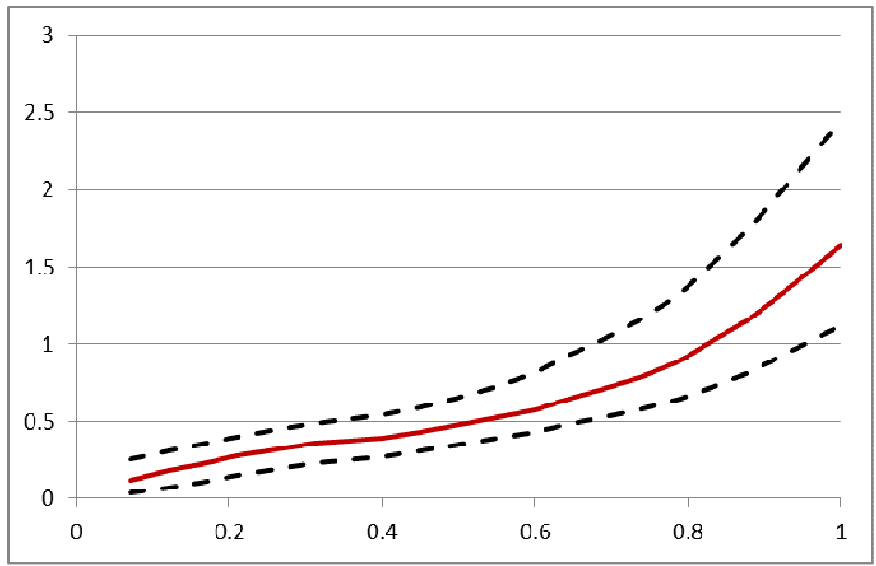

**Figure 1:** MML-interval for the quantiles of the STS distribution
($n = 20$, $\alpha = 0.05$, $\mu = 0, \sigma = 1, r = 2, d = 0$)

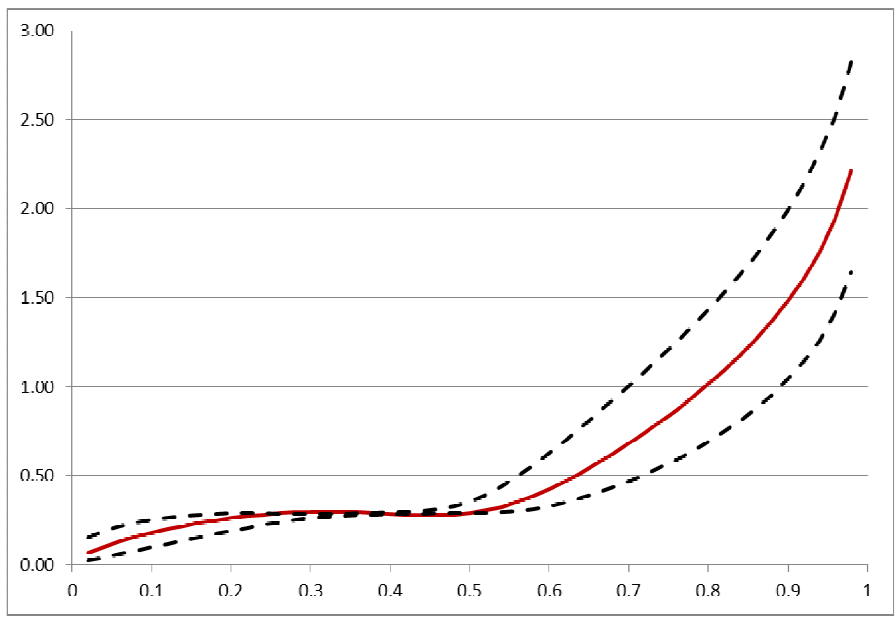

**Figure 2:** MML-interval for the quantiles of the STS distribution
($n = 50$, $\alpha = 0.05$, $\mu = 0, \sigma = 1, r = 2, d = 1.0$)

## 4.2. Long-Tailed Symmetric Distribution

The pdf for the long-tailed symmetric (LTS) distribution is given by

$$f(y,p) \propto \frac{1}{\sigma}\left\{1 + \frac{(y-\mu)^2}{k\sigma^2}\right\}^{-p}, \quad -\infty < y < \infty; \tag{22}$$

with mean $\mu$ and variance $\sigma^2$; $k = 2p-3$ and $p \geq 2$ (see Tiku and Suresh, 1992). The MML estimators of a complete sample of size $n$ are obtained as explained in Section 2 and they have the same structures given in (15) and (16) with

$$B = \frac{2p}{k}\sum_{i=1}^{n}\alpha_{1i}y_{(i)}, \quad C = \frac{2p}{k}\left\{\sum_{i=1}^{n}\beta_{1i}y_{(i)}^2 - m\hat{\mu}^2\right\} \tag{23}$$

where

$$\alpha_{1i} = \frac{(2/k)t_i^2}{[1+(1/k)t_i^2]^2}, \quad \beta_{1i} = \frac{1-(1/k)t_i^2}{[1+(1/k)t_i^2]^2} \quad \text{and} \quad m = \frac{2p\sum_{i=1}^{n}\beta_{1i}}{k}. \tag{24}$$

Expected values of the standardized order statistics ($t_i$) are available in Tiku and Kumra (1981) for p=2(0.5)10, n≤20. One can also see the corresponding expected values for $p=1.5$, n≤20 in Vaughan (1992b) and those for $p=1$ (Cauchy distribution), (n≥6) in Vaughan (1994). Tiku et. al. (2000) provides the alternative coefficients $\alpha_{2i}$ and $\beta_{2i}$ in case $\beta_{1i} < 0$ for some $i$ and they are given by

$$\alpha_i^* = 0 \quad \text{and} \quad \beta_i^* = 1/\{1+(1/k)t_{(i)}^2\}. \tag{25}$$

Table 3 and 4 show the results of a simulation study for the efficiencies of the hazard rate estimators obtained for the LTS distribution. As in the case of the STS distribution, the same conclusion can be drawn and the estimators based on the MML estimation perform much more efficient than those based on the LS. As for the coverage probabilities, results are similar to those obtained from the STS distribution.

**Table 3:** Means and Variances of Hazard Rate Estimators for the Quantiles of the LTS Distribution (Q) ($n = 20$, $\mu = 0, \sigma = 1, p = 3$)

$\hat{\mu} = 0.0019$, $\hat{\sigma} = 1.049$, $\bar{x} = 0.0018$, $s = 0.972$

| | | Means for n=20 | | | | Variances for n=20 | | | |
|---|---|---|---|---|---|---|---|---|---|
| Q | Exact | $HR^1_{MML}$ | $HR^1_{LS}$ | $HR^2_{MML}$ | $HR^2_{LS}$ | $HR^1_{MML}$ | $HR^1_{LS}$ | $HR^2_{MML}$ | $HR^2_{LS}$ |
| 0.07 | 0.0532 | 0.0628 | 0.0538 | 0.0584 | 0.0467 | 0.0012 | 0.0014 | 0.0014 | 0.0017 |
| 0.12 | 0.0948 | 0.1039 | 0.0922 | 0.1008 | 0.0876 | 0.0020 | 0.0024 | 0.0022 | 0.0026 |
| 0.17 | 0.1338 | 0.1406 | 0.1282 | 0.1385 | 0.1252 | 0.0026 | 0.0031 | 0.0027 | 0.0033 |
| 0.21 | 0.1701 | 0.1740 | 0.1620 | 0.1728 | 0.1606 | 0.0027 | 0.0034 | 0.0029 | 0.0037 |
| 0.26 | 0.2032 | 0.2041 | 0.1934 | 0.2038 | 0.1935 | 0.0026 | 0.0034 | 0.0028 | 0.0037 |
| 0.31 | 0.2329 | 0.2308 | 0.2218 | 0.2312 | 0.2232 | 0.0023 | 0.0031 | 0.0025 | 0.0033 |
| 0.36 | 0.2588 | 0.2540 | 0.2468 | 0.2551 | 0.2495 | 0.0019 | 0.0025 | 0.0020 | 0.0027 |
| 0.40 | 0.2806 | 0.2734 | 0.2679 | 0.2749 | 0.2715 | 0.0014 | 0.0019 | 0.0015 | 0.0021 |
| 0.45 | 0.2978 | 0.2887 | 0.2845 | 0.2906 | 0.2890 | 0.0009 | 0.0013 | 0.0010 | 0.0014 |
| 0.50 | 0.3102 | 0.2996 | 0.2961 | 0.3019 | 0.3013 | 0.0005 | 0.0008 | 0.0006 | 0.0009 |
| 0.55 | 0.3174 | 0.3058 | 0.3024 | 0.3083 | 0.3079 | 0.0003 | 0.0005 | 0.0004 | 0.0006 |
| 0.60 | 0.3188 | 0.3069 | 0.3028 | 0.3094 | 0.3083 | 0.0003 | 0.0005 | 0.0003 | 0.0006 |
| 0.64 | 0.3139 | 0.3026 | 0.2969 | 0.3050 | 0.3024 | 0.0005 | 0.0008 | 0.0006 | 0.0009 |
| 0.69 | 0.3023 | 0.2923 | 0.2844 | 0.2945 | 0.2896 | 0.0010 | 0.0015 | 0.0011 | 0.0016 |
| 0.74 | 0.2830 | 0.2754 | 0.2648 | 0.2770 | 0.2691 | 0.0016 | 0.0023 | 0.0018 | 0.0026 |
| 0.79 | 0.2553 | 0.2512 | 0.2378 | 0.2520 | 0.2407 | 0.0024 | 0.0032 | 0.0026 | 0.0036 |
| 0.83 | 0.2179 | 0.2185 | 0.2026 | 0.2183 | 0.2038 | 0.0031 | 0.0040 | 0.0033 | 0.0044 |
| 0.88 | 0.1692 | 0.1754 | 0.1582 | 0.1738 | 0.1567 | 0.0034 | 0.0041 | 0.0037 | 0.0046 |
| 0.93 | 0.2728 | 0.2665 | 0.2547 | 0.2678 | 0.2585 | 0.0020 | 0.0027 | 0.0021 | 0.0030 |

**Table 4:** Coverage Probabilities $(1-\alpha)$ of the CI's for the Hazard Rate Estimators for the LTS Distribution ($n = 20$, $\mu = 0, \sigma = 1, p = 3$)

$\hat{\mu} = 0.0019$, $\hat{\sigma} = 1.049$, $\bar{x} = 0.0018$, $s = 0.972$

| | | Coverage Probabilities | | | |
|---|---|---|---|---|---|
| Q | $HR^2_{MML}$ | $HR^2_{LS}$ | Q | $HR^2_{MML}$ | $HR^2_{LS}$ |
| 0.07 | 0.771 | 0.741 | 0.55 | 0.955 | 0.953 |
| 0.12 | 0.840 | 0.817 | 0.59 | 0.952 | 0.949 |
| 0.17 | 0.880 | 0.858 | 0.64 | 0.949 | 0.940 |
| 0.21 | 0.906 | 0.888 | 0.69 | 0.939 | 0.927 |
| 0.26 | 0.922 | 0.909 | 0.74 | 0.926 | 0.912 |
| 0.31 | 0.933 | 0.923 | 0.79 | 0.906 | 0.888 |
| 0.36 | 0.943 | 0.937 | 0.83 | 0.880 | 0.857 |
| 0.40 | 0.950 | 0.947 | 0.88 | 0.839 | 0.814 |
| 0.45 | 0.953 | 0.953 | 0.93 | 0.918 | 0.904 |
| 0.50 | 0.956 | 0.955 | | | |

### 4.3. Other Location-Scale Distributions

Similarly to the STS and the LTS distributions, we can easily derive the hazard rate estimators for other location-scale distributions having monotone hazard rates and their CI's. MML estimators of unknown parameters for these distributions are not difficult to obtain due to the efficient approximation of nonlinear terms of likelihood functions. Once MML estimators are obtained, it is again very easy to estimate the hazard rate at specific point $t$. For location-scale distributions having nonmonotonic hazard rates, we can still make use of the same approach. Except at the change point, it is possible to obtain highly accurate estimates of hazard rates. As a future study, the change point problem can be worked out in more detail and new approximation methods can be developed to increase accuracies.

### 5. Estimation in Censored Samples

For censored samples coming from location-scale families, unknown location and scale estimators can easily be obtained through the MML approach as explained in Section 2. For the hazard rate estimation, however, the procedure proposed in Section 3 is not directly applicable since the order statistics in the censored parts of samples cannot be observed. For a future $t$ point, we cannot observe two order statistics that will be utilized in the linear approximation. This leads us to propose two new approaches which will be helpful in predicting order statistics for the linearization boundaries.

Consider a censored sample from a location-scale family and assume last $n-r$ observations are not observed. In order to linearize the hazard function at point $t$ $(t > x_{(r)})$, we need to have an interval of two consecutive order statistics containing $t$. Since these order statistics are not available, we can proceed with their predicted values. To do that, consider the following predictive log-likelihood function of the censored sample and the $i$th order statistic $x_{(i)}$ $(i > r)$:

$$\ln L(X|\mu,\sigma,x_{(i)}) \propto \sum_{i=1}^{r} \ln f\left(\frac{x_{(j)}-\mu}{\sigma}\right) + \ln f\left(\frac{x_{(i)}-\mu}{\sigma}\right) + (i-r-1)\ln\left[F\left(\frac{x_{(i)}-\mu}{\sigma}\right) - F\left(\frac{x_{(r)}-\mu}{\sigma}\right)\right]$$
$$+ (n-i)\ln\left[1-F\left(\frac{x_{(i)}-\mu}{\sigma}\right)\right]; \quad r < i \leq n \quad (26)$$

Taking the partial derivatives with respect to $\mu, \sigma$ and $x_{(i)}$ gives nonlinear terms as was discussed in Section 2. However, the MML approach easily takes care of the problem and provides closed-form solutions of the estimators $\hat{\mu}$, $\hat{\sigma}$ and $\hat{x}_{(i)}$. See also Raqab (1997) who also considered the same problem by using the similar approximation techniques. Once the future order statistics are predicted, hazard rate estimation at point $t$ can easily be handled by simply determining the value of $k$ satisfying $\hat{x}_{(k)} < t < \hat{x}_{(k+1)}$. Then, the estimation follows from the procedure explained in Section 3.

Another way of determining the value of $k$ is to make use of the generalized spacings of the location-scale distribution since the spacings are unbiased estimators of the scale parameter $\sigma$. Let $\hat{\mu}$ and $\hat{\sigma}$ be the MML estimators obtained by the MML approach for the censored sample $x_{(1)}, x_{(2)}, x_{(3)}, ..., x_{(r)}$. Then, an estimator of a future statistic $x_{(i)}$ is

$$\hat{x}_{(i)} = \frac{\sum_{j=1}^{r}(x_{(j)} + \hat{\sigma}(t_i - t_j))}{r} \tag{27}$$

where $t_i$ is the expected value of the $i$th standardized order statistic. One can work on the efficiency of this estimator as a future study.

**References**

Akkaya, A. D. and Tiku, M. L., (2008). Short-tailed Distributions and Inliers, Test, 17, 282-299.

Antoniadis, A., Gregoire, G. (1999), "Density and Hazard rate estimation for right-censored data using wavelet methods," *J. R. Statist. Soc. B*, 61, 1, pp. 63–84.

Bezandry, P. H., Bonney, G. E. and Gannoun, A. (2005). Consistent estimation of the density and hazard rate functions for censored data via the wavelet method. Statistics and Probability Letters, 74(4), 366-372.

Brunel E, Comte F (2008) Adaptive estimation of hazard rate with censored data. Comm Statist Theory Methods 37(8–10):1284–1305.

Mahapatra, A. K., Kumar, S. and Vellaisamy, P. (2012). Simultaneous estimation of hazard rates of several exponential populations, Statistica Neerlendica, 66(2), 121-132.

Meeker, W. Q., Escobar L. A., and Hill, D. A. (1992). Sample sizes for estimating the Weibull hazard function from censored samples, IEEE Transactions on Reliability R-41, 133-138.

Rai, B. And Singh, N. (2003). Hazard rate estimation from incomplete and unclean warranty data, Rel. Eng. Sys. Saf., 81, 79-92.

Raqab, M. Z. (1997). Modified Maximum Likelihood Predictors of Future Order Statistics from Normal Samples, *Computational Statistics and Data Analysis* 25, 91-106.

Tiku, M.L. (1967). Estimating the mean and standard deviation from a censored normal sample. Biometrika, **54**, 155-165.

Tiku, ML. and Akkaya, A.D. (2004). Robust Estimation and Hypothesis Testing, New Age International Publishers (P): New Delhi, 2004.

Tiku, M.L. and Kumra, S. (1981) Expected values and variances and covariances of order statistics for a family of symmetrical distributions (Student's t). Selected Tables in Mathematical Statistics 8, American Mathematical Society, Providence, RI: 141-270.

Tiku, M. L. and Suresh, R. P. (1992). A new method of estimation for location and scale parameters. J. Stat. Plann. Inf. 30, 281-292.

Tiku, M.L. and Vaughan, D.C. (1999). A family of short-tailed symmetric distributions. Technical Report: McMaster University, Canada.

Tiku, M.L., Wong, W.K., Vaughan, D.C. and Bian, G.R. (2000). Time series models in non-normal situations: Symmetric innovations", J. Time Series Analysis 21, 451-469.

Vaughan, D.C. (1992b) Expected values, variances and covariances of order statistics for Student's t distribution with two degrees of freedom. Commun. Stat. Simul. 21, 391-404.

Vaughan, D.C. (1992). On the Tiku-Suresh method of estimation. Commun. Statist.-Theory Meth., **21**, 451-469.

Vaughan, D.C. (1994) The exact values of the expected values, variances and covariances of the order statistics from the Cauchy distribution. J. Stat. Comput. Simul. 49, 21-32.